\documentclass[12pt]{article}
\usepackage{pifont}

\usepackage{amsfonts}
\usepackage{latexsym}
\usepackage{amsmath}
\usepackage{amssymb}
\usepackage{color}

\def\d{\delta}
\def\a{\alpha}
\def\p{\varphi}

\def\e{\varepsilon}
\def\z{\zeta}
\def\P{\Phi}
\def\b{\beta}

\title{The optimal control related to Riemannian manifolds and the viscosity solutions to H-J-B equations}

\date{}

 \author{\  Xuehong Zhu \\
\small{Institute of Mathematics, \ Shandong University }\\
\small{ Jinan, \ 250100, \ China}\\
\small {School of Science, \ Nanjing}\\
 \small{University of
Aeronautics and
Astronautics }\\
\small{ Nanjing, \ 210016, \ China}\\
\small{{E-mail: hilda2002@163.com }}}

\begin{document}

\maketitle

\begin{abstract}
This paper is concerned with the Dynamic Programming Principle (DPP
in short)
  with SDEs on Riemannian manifolds. Moreover, through the DPP, we
  conclude that the cost function is the unique viscosity solution to the
  related PDEs on manifolds.\\
\par $\textit{Keywords:}$ Dynamic programming principle; Riemannian manifold;
Viscosity solution.
\end{abstract}



\section{Introduction}\label{sec:intro}

El Karoui, Peng and Quenez [3] gave the formulation of recursive
utilities and their properties from the BSDE point of view. As we
know, the recursive optimal control problem is represented as a kind
of optimal control problem whose cost functional is described by the
solution of BSDE. In 1992, Peng [5] got the Bellman's dynamic
programming principle for this kind of problem and proved that the
value function is a viscosity solution of one kind of quasi-linear
second-order partial differential equation (PDE in short) which is
the well-known as Hamilton-Jacobi-Bellman (H-J-B in short) equation.
Later in 1997, he virtually generalized these results to a much more
general situation, under Markvian and even Non-Markvian framework
([6], Chapter 2).

But sometimes, in financial market, as the solution to a SDE with
control, the wealth process of the investor may be constrained, for
example, it should be nonnegative. In particular, for some special
need, it may be a process in some curving spaces. So it is natural
to consider the following question: if the SDE in stochastic
recursive optimal control problems is defined on Riemannian
manifolds, then will we still have the similar results as what we
have mentioned in $R^{n}$? The objective of this paper is to study
this problem.

\vspace{3mm}

Let $(W(t),t\geq 0)$ be a $d$-dimensional standard Brownian motion
  on some complete probability space $(\Omega,{\cal{F}},P)$. We
  denote by $({\cal{F}}_{t})_{t\geq 0}$ the natural filtration
  generated by $W$ and augmented by the $P$-null sets of
  ${\cal{F}}$.

Let $U$ be a compact subset of $R^{d+1}$. We call a function $f :
\Omega\times [t,T]\rightarrow U$ an admissible control if it's an
adapted stochastic process. We denote by  $\mathcal{U}_{t,T}$ the
set of all admissible controls.

Assume that $M$ is a compact Riemmannian manifold without boundary.
Now we can consider the following controlled stochastic differential
equation on $M$ in a fixed time interval $[t,T]$:
$$
\left\{\aligned &
dX_{s}^{t,\z;v.}=v_{0}(s)V_{0}(s,X_{s}^{t,\z;v.})ds+\sum_{\alpha=1}^{d}V_{\a}(s,X_{s}^{t,\z;v.})\circ
v_{\alpha}(s)dW^{\a}_{s},\cr\noalign{\vskip 2mm} &
X_{t}^{t,\z;v.}=\z\in M,\endaligned\right. \eqno(1.1)
$$
where $\z$ is $\mathcal{F}_{t}$-measurable,
$v.=v(\cdot):=(v_{0}(\cdot),v_{1}(\cdot),...,v_{d}(\cdot))\in
\mathcal{U}_{t,T}$, and $V_{0},V_{1},...V_{d}$ are $d+1$
deterministic one-parameter smooth vector fields on $M$.

Since $M$ is compact and without boundary, according to [4], there
exists a unique $M$-valued continuous process which solves equation
(1.1). Moreover, this solution does not explode.

Let us consider functions $f : [0,T]\times M\times R\times
R^{1\times d}\times U\rightarrow R$ and $\Phi : M\rightarrow R$
which satisfy:

(A1).there exists a constant $K\geq 0$, s.t., we have: $\forall t,
\forall (x,y,z,v)\mbox{ \ and \ }(x',y',z',v'),$
$$|\P(x)-\P(x')|+|f(t,x,y,z,v)-f(t,x',y',z',v')|\leq
K(|y-y'|+|z-z'|+d(x,x')+|v-v'|),
$$

(A2).there exists a constant $K_{0}\geq 0$, s.t.,  $\forall
(t,x,v), |f(t,x,0,0,v)|\leq K_{0}$, \\
where $d(\cdot,\cdot)$ denotes the Riemannian distance function on
$M$.

By the above assumptions, according to [6], there exists a unique
solution $(Y., Z.)\in \mathcal{M}(t,T; R\times R^{1\times d})$ to
the following BSDE:
$$
\left\{\begin{array}{l} -dY_{s}^{t,\zeta ;v.}=f(s,X_{s}^{t,\zeta
;v.},Y_{s}^{t,\zeta ;v.},Z_{s}^{t,\zeta ;v.},v_{s})-Z_{s}^{t,\zeta
;v.}dW_{s},s\in \lbrack t,T],
\\
Y_{T}^{t,\zeta ;v.}=\Phi (X_{T}^{t,\zeta ;v.}),%
\end{array}\right.
$$
where $\mathcal{M}(0,T;R^{n})$ denotes the Hilbert space of adapted
stochastic processes $f :\Omega\times [0,T]\rightarrow R^{n}$ such
that
$$
\|f\|=( E\int_{0}^{T}|f(t)|^{2}dt)^{\frac{1}{2}}<\infty.
$$

 When $\z=x\in M$ is deterministic, We define
$$
J(t,x;v(\cdot)):=Y_{s}^{t,x ;v.}|_{s=t}.
$$
This is the so-called cost function. And then we can define a value
function of the optimal control problem as follows:
$$
u(t,x):=\mbox{essinf}_{v(\cdot)\in
\mathcal{U}_{t,T}}J(t,x;v(\cdot)).$$

Our purpose is to get the general Dynamic Programming Principle of
the value function $u(t,x)$.

 \vskip 2mm

\section{Dynamic Programming Principle}\label{sec:intro}

\vskip 2mm

If we define
$$\mathcal{U}_{t,T}^{t}:=\{v(\cdot)\in \mathcal{U}_{t,T} : v(\cdot) \mbox{ \ is \ }
\mathcal{F}^{t}_{s}-\mbox{adapted}\},$$ where
$\mathcal{F}^{t}_{s}:=\sigma\{W_{r}-W_{t} , t\leq r\leq s\}.$

By Proposition 5.1 in [6], there exist
$\{v^{i}(\cdot)\}^{\infty}_{i=1},v^{i}(\cdot)\in
\mathcal{U}_{t,T}^{t}$, such that
$u(t,x)=\lim_{i\rightarrow\infty}J(t,x;v^{i}(\cdot))$ and $u(t,x)$
is a deterministic function, i.e.,
$$
u(t,x):=\mbox{essinf}_{v(\cdot)\in
\mathcal{U}_{t,T}}J(t,x;v(\cdot))=\inf_{v(\cdot)\in
\mathcal{U}_{t,T}^{t}}J(t,x;v(\cdot)).
$$

Since our SDE is defined on Riemannian manifolds, $\forall
(\z,v(\cdot))$ and $(\z',v'(\cdot))$,\\
$d^{2}(X^{t,\z;v.}_{s},X^{t,\z';v'.}_{s})$ is not necessarily twice
differentiable. So the good estimate about the continuous dependence
of $X.^{t,\z;v.}$ w.r.t to $(\z,v.)$ does not hold and nor does
$Y_{s}^{t,\z;v.}|_{s=t}$. They're unfavourable factors for our
dynamic programming principle. So we turn to the embedding mapping.

By the Whitney's theorem, there exists an embedding mapping $\Psi$
such that, $\Psi : M\rightarrow \Psi(M)\subset R^{n}$ for some $n\in
N$. Set $\Psi(X.^{t,\z;v.})=\tilde{X}.^{t,\tilde{\z};v.}$, where
$\tilde{\z}=\Psi(\z)$. Then $\tilde{X}.^{t,\tilde{\z};v.}$ satisfies
the following SDE on $\Psi(M)$:
$$
\left\{\aligned &
d\tilde{X}_{s}^{t,\tilde{\z};v.}=v_{0}(s)\tilde{V}_{0}(s,\tilde{X}_{s}^{t,\tilde{\z};v.})ds+\sum_{\alpha=1}^{d}\tilde{V}_{\a}(s,\tilde{X}_{s}^{t,\tilde{\z};v.})\circ
v_{\alpha}(s)dW^{\a}_{s},\cr\noalign{\vskip 2mm} &
\tilde{X}_{t}^{t,\tilde{\z};v.}=\tilde{\z}\in
\Psi(M),\endaligned\right. \eqno(2.1)
$$
where $\tilde{V}_{\a}=\Psi_{*}V_{\a}, \a=0,1,...,d$ and $\Psi_{*}$
is the tangent mapping. And we can extend each $\tilde{V}_{\a}$ to
smooth vector field defined on $R^{n}$ with compact support. We
denote the extensions still by $\tilde{V}_{\a}$.

So we have the following SDE in $R^{n}$ whose initial point is in
$\Psi(M)$:
$$
\left\{\aligned &
d\tilde{X}_{s}^{t,\tilde{\z};v.}=v_{0}(s)\tilde{V}_{0}(s,\tilde{X}_{s}^{t,\tilde{\z};v.})ds+\sum_{\alpha=1}^{d}v_{\alpha}(s)\tilde{V}_{\a}(s,\tilde{X}_{s}^{t,\tilde{\z};v.})
dW^{\a}_{s}+\frac{1}{2}\sum_{\alpha=1}^{d}v^{2}_{\a}(s)\nabla_{\tilde{V_{\a}}}\tilde{V_{\a}}(s,\tilde{X}_{s}^{t,\tilde{\z};v.})ds,\cr\noalign{\vskip
2mm} & \tilde{X}_{t}^{t,\tilde{\z};v.}=\tilde{\z}\in
\Psi(M),\endaligned\right. \eqno(2.2)
$$
where $\nabla$ is the connection of $R^{n}$. According to [4], SDE
(2.2) has the same unique solution with SDE (2.1), i.e., although
although SDE (2.2) is defined in $R^{n}$, as long as its initial
point is in $\Psi(M)$, it won't leave $\Psi(M)$.

Since $U$ is compact and each $\tilde{V}_{i}$ is a smooth vector
field in $R^{n}$ with compact support, the coefficients of SDE (2.2)
are bounded and Lipschitz continuous w.r.t. $x$ and $v$. By [6], we
have the following estimate:
$$
E^{\mathcal{F}_{t}}[\sup_{s\in[t,T]}|\tilde{X}_{s}^{t,\tilde{\z};v.}-\tilde{X}_{s}^{t,\tilde{\z'};v'.}|^{2}]\leq
C(|\tilde{\z}-\tilde{\z'}|^{2}+E^{\mathcal{F}_{t}}\int_{t}^{T}|v(s)-v'(s)|^{2}ds),
\eqno(2.3)
$$
where $C>0$ is a constant which only depends on the Lipschitz
constant of the coefficients of SDE (2.2). Here and in the sequel,
the constant $C$ appearing in each estimate won't be necessarily the
same one.

{\bf Lemma 2.1. }$Y_{t}^{t,x;v.}$ is continuous w.r.t $(x,
v(\cdot))$ and it is uniformly continuous in $x$, uniformly in
$(t,v(\cdot))$.

{\bf Proof: }Using It\^{o}'s formula to
$e^{\b(s-t)}|Y^{t,x;v.}_{s}-Y^{t,x';v'.}_{s}|^{2}$ for some positive
constant $\b$, we have
$$
\begin{array}{ll}
&
|Y^{t,x;v.}_{t}-Y^{t,x';v'_.}_{t}|^{2}+E^{\mathcal{F}_{t}}\displaystyle\int_{t}^{T}e^{\b(s-t)}[\b
|Y^{t,x;v.}_{s}-Y^{t,x';v'_.}_{s}|^{2} +
|Z^{t,x;v.}_{s}-Z^{t,x';v'_.}_{s}|^{2}]ds\\
=&
e^{\b(T-t)}E^{\mathcal{F}_{t}}|\P(X^{t,x;v.}_{T})-\P(X^{t,x';v'_.}_{T})|^{2}
+E^{\mathcal{F}_{t}}\displaystyle\int_{t}^{T}
e^{\b(s-t)}2(Y^{t,x;v.}_{s}-Y^{t,x';v'_.}_{s})*\\
&\ \ \ \ \ \ \ \ \ \ \ \ \ \ \ \ (f(s,X^{t,x;v.}_{s},Y^{t,x;v.}_{s},Z^{t,x;v.}_{s},v_{s})-f(s,X^{t,x';v'_.}_{s},Y^{t,x';v'_.}_{s},Z^{t,x';v'_.}_{s},v'_{s}))ds\\
\leq &
Ke^{\b(T-t)}E^{\mathcal{F}_{t}}d^{2}(X^{t,x;v.}_{T},X^{t,x';v'_.}_{T})+\frac{1}{2}E^{\mathcal{F}_{t}}\int_{t}^{T}
e^{\b(s-t)}|Z^{t,x;v.}_{s}-Z^{t,x';v'_.}_{s}|^{2}ds\\
&+E^{\mathcal{F}_{t}}\displaystyle\int_{t}^{T}
e^{\b(s-t)}d^{2}(X^{t,x;v.}_{s},X^{t,x';v'_.}_{s})ds+E^{\mathcal{F}_{t}}\displaystyle\int_{t}^{T}
e^{\b(s-t)}|v_{s}-v'_{s}|^{2}ds\\
&+E^{\mathcal{F}_{t}}\displaystyle\int_{t}^{T}
e^{\b(s-t)}(2K+2K^{2}+K^{2}+K^{2})|Y^{t,x;v.}_{s}-Y^{t,x';v'_.}_{s}|^{2}ds.
\end{array}
$$
If we choose $\b=2K+4K^{2}+1$, we have
$$
\begin{array}{ll}
&
|Y^{t,x;v.}_{t}-Y^{t,x';v'_.}_{t}|^{2}+E^{\mathcal{F}_{t}}\displaystyle\int_{t}^{T}
e^{\b(s-t)}[|Y^{t,x;v.}_{s}-Y^{t,x';v'_.}_{s}|^{2}+\frac{1}{2}|Z^{t,x;v.}_{s}-Z^{t,x';v'_.}_{s}|^{2}]ds\\
\leq & C\{E^{\mathcal{F}_{t}}[\sup_{s\in
[t,T]}d^{2}(X^{t,x;v.}_{s},X^{t,x';v'_.}_{s})]\\
& \ \ \ \ +E^{\mathcal{F}_{t}}\displaystyle\int_{t}^{T}\sup_{s\in
[t,T]}d^{2}(X^{t,x;v.}_{s},X^{t,x';v'_.}_{s})ds+E^{\mathcal{F}_{t}}\displaystyle\int_{t}^{T}|v_{s}-v'_{s}|^{2}ds\}.
\end{array}\eqno{(2.4)}
$$

Since $M$ is compact, $\Psi : M\rightarrow \Psi(M)$ and $\Psi^{-1} :
\Psi(M)\rightarrow M$ are both uniformly continuous mappings. So
with (2.3), when $(x',v')\rightarrow (x,v)$, we have $$
E^{\mathcal{F}_{t}}[\sup_{s\in[t,T]}d^{2}(X^{t,x;v.}_{s},X^{t,x';v'.}_{s})]\rightarrow
0 .$$ Moreover, when $x'\rightarrow x$,
$$
\sup_{v(\cdot)\in
\mathcal{U}_{t,T}}E^{\mathcal{F}_{t}}[\sup_{s\in[t,T]}d^{2}(X^{t,x;v.}_{s},X^{t,x';v.}_{s})]\leq
E^{\sup_{v(\cdot)\in
\mathcal{U}_{t,T}}\mathcal{F}_{t}}[\sup_{s\in[t,T]}d^{2}(X^{t,x;v.}_{s},X^{t,x';v.}_{s})]
\rightarrow 0.\eqno{(2.5)}
$$
Through the theorem of control convergence, we have that
$Y_{t}^{t,x;v.}$ is continuous w.r.t $(x, v(\cdot))$.

What's more, $\forall \e>0$, we choose $\d_{0}=\frac{\e}{C(1+T)}$.
By (2.5), for this $\d_{0}$, there exists $\d>0$, such that, when
$d(x,x')<\d$ (here $\d$ doesn't depend on $x$ or
$x'$),$$\sup_{v(\cdot)\in
\mathcal{U}_{t,T}}E^{\mathcal{F}_{t}}[\sup_{s\in[t,T]}d^{2}(X^{t,x;v.}_{s},X^{t,x';v.}_{s})]<\d_{0}.$$
Combining (2.4), we have
$$
\begin{array}{ll}
&|Y^{t,x;v.}_{t}-Y^{t,x';v.}_{t}|^{2}+E^{\mathcal{F}_{t}}\int_{t}^{T}
e^{\b(s-t)}[|Y^{t,x;v.}_{s}-Y^{t,x';v.}_{s}|^{2}+\frac{1}{2}|Z^{t,x;v.}_{s}-Z^{t,x';v.}_{s}|^{2}]ds\\
<& C\d_{0}+CT\d_{0}=\e.
\end{array}
$$
So we have finished the proof.
$$
\eqno{\Box}
$$

And we can get some properties of $u(t,x)$:

{\bf Lemma 2.2. }$u(t,x)$ is bounded and uniformly continuous in
$x$, uniformly in $t$.

{\bf Proof: }Applying It\^{o}'s formula to
$e^{\b_{1}(s-t)}|Y^{t,x;v.}_{s}|^{2}(\b_{1}=2K+2K^{2}+2)$, with the
same method of above, we have
$$
\begin{array}{ll}
&|Y^{t,x;v.}_{t}|^{2}+E^{\mathcal{F}_{t}}\int_{t}^{T}\b_{1}
e^{\b_{1}(s-t)}|Y^{t,x;v.}_{s}|^{2}ds+E^{\mathcal{F}_{t}}\displaystyle\int_{t}^{T}
e^{\b_{1}(s-t)}|Z^{t,x;v.}|^{2}ds\\
\leq & e^{\b_{1}(T-t)}E^{\mathcal{F}_{t}}|\P(X^{t,x;v.}_{T})|\\
&+ E^{\mathcal{F}_{t}}\displaystyle\int_{t}^{T}
e^{\b_{1}(s-t)}2|Y^{t,x;v.}_{s}|(K|Y^{t,x;v.}_{s}|+K|Z^{t,x;v.}_{s}|+|f(s,X^{t,x;v.}_{s},0,0,v_{s})|)ds.
\end{array}
$$
Since $\P(\cdot)$ is continuous and $M$ is compact, $\P(x)$ is
bounded. This with A2, we have that there exists some constant $C$
independent of $(t,x,v(\cdot))$, such that,
$$
|Y^{t,x;v.}_{t}|^{2}+E^{\mathcal{F}_{t}}\int_{t}^{T}
e^{\b_{1}(s-t)}|Y^{t,x;v.}_{s}|^{2}ds+\frac{1}{2}E^{\mathcal{F}_{t}}\int_{t}^{T}
e^{\b_{1}(s-t)}|Z^{t,x;v.}|^{2}ds \leq C. \eqno{(2.6)}
$$
So $u(t,x)\leq C$.

By the definition of $u(t,x)$, we know that for any $\e>0$, there
exist $v(\cdot), v'(\cdot)\in \mathcal{U}_{t,T}^{t}$, such that,
$$Y^{t,x;v.}_{t}-\e\leq u(t,x),\mbox{ \ }Y^{t,x';v'.}_{t}-\e\leq u(t,x').$$
So we have
$$
Y^{t,x;v.}_{t}-\e\leq u(t,x)\leq Y^{t,x;v'.}_{t},\mbox{ \
}Y^{t,x';v'.}_{t}-\e\leq u(t,x')\leq Y^{t,x';v.}_{t}.
$$
And that yields
$$
Y^{t,x;v.}_{t}-Y^{t,x';v.}_{t}-\e\leq u(t,x)-u(t,x')\leq
Y^{t,x;v'.}_{t}-Y^{t,x';v'.}_{t}+\e.
$$
Since $Y^{t,x;v.}_{t}$ is continuous in $x$ uniformly in
$(t,v(\cdot))$, we've got our conclusion.
$$
\eqno{\Box}
$$

If we replace the variable $x$ in $u(t,x)$ by a r.v. $\z$ which is
$\mathcal{F}_{t}$-measurable, we have:

{\bf Lemma 2.3. }For any fixed $t\in [0,T]$ and $\z$ which is
$\mathcal{F}_{t}$-measurable, we have:

(i)$\forall v(\cdot)\in \mathcal{U}_{t,T}, u(t,\z)\leq
Y^{t,\z;v.}_{t},$

(ii)$\forall \e>0,$ there exists a $v(\cdot) \in \mathcal{U}_{t,T}$
such that $u(t,\z)\geq Y^{t,\z;v.}_{t} -\e.$

{\bf Proof: }We have known that $u$ is continuous in $x$ and
$Y^{t,\z;v.}_{t}$ is continuous in $(\z,v(\cdot))$. Recall that the
collection of processes $(v(s))_{s\in [t,T]}$ with
$$
\{v(s)=\sum_{i=1}^{N}I_{A_{i}}v^{i}(s) : \{A_{i}\}_{i=1}^{N} \mbox{
\ is a \ } \mathcal{F}_{t}-\mbox{partition of \ }\Omega,
v^{i}(\cdot)\in \mathcal{U}_{t,T} \mbox{ \ is \
}\mathcal{F}^{t}_{s}-\mbox{adapted}\}
$$
is dense in $\mathcal{U}_{t,T}$. So for (i), we need only to discuss
special $\z$ and $v(\cdot)$ as follows:
$$
\z=\sum_{i=1}^{N}I_{A_{i}}x_{i},
v(\cdot)=\sum_{i=1}^{N}I_{A_{i}}v^{i}(\cdot),
$$
where $A_{i}$ and $v^{i}(\cdot)$ are described as above, and
$x_{i}\in M,i=1,...,N$. Then we can use the same method with Theorem
4.7 in [6] to get
$$
Y^{t,\z;v.}_{t}=\sum_{i=1}^{N}I_{A_{i}}Y^{t,x_{i};v^{i}.}_{t}\geq
\sum_{i=1}^{N}I_{A_{i}}u(t,x_{i})=u(t,\sum_{i=1}^{N}I_{A_{i}}x_{i})=u(t,\z).
$$

For (ii), we can use the same technique. Considering that $M$ is
compact, for any $n\in N$, there exist a collection of data
$\{U_{i},\p_{i}\}_{i=1}^{N_{n}}$ such that $\dim(U_{i})<
\frac{1}{2^{n}}$, where $\dim(U_{i}):=\sup_{x,y\in M}d(x,y)$.  For
any $\z$ which is $\mathcal{F}_{t}$-measurable, choose any fixed
$y_{i}\in U_{i}$ and set
$\eta_{n}=\sum_{i=1}^{N_{n}}y_{i}I_{\{\omega : \z(\omega)\in
U_{i}\}}$. Then we have
$$
d(\eta_{n},\z)<\frac{1}{2^{n}},P-a.s.
$$
By Lemma 2.1. and 2.2., for any $\e>0$, there exists $\d>0$, such
that, when $d(x,x')<\d$,
$$
u(t,x)-u(t,x')\geq -\frac{\e}{3}, \mbox{ \ }
Y^{t,x;v}_{t}-Y^{t,x';v}_{t}\geq -\frac{\e}{3}, \mbox{ \ for any \
}(t,v(\cdot))\in [0,T]\times \mathcal{U}_{t,T}.
$$
For this $\d$, through the discussion above, there exists
$\eta=\sum_{i=1}^{N_{\d}}x_{i}I_{A_{i}} \mbox{ \ where \ } x^{i}\in
M \mbox{ \ and \ } \{A_{i}\}_{i=1}^{N_{\d}} \mbox{ \ is a \ }
\mathcal{F}_{t}-\mbox{partition of \ }\Omega,$ such that
$d(\eta,\z)<\d,P-a.s.$. So we have, for any $(t,v(\cdot))\in
[0,T]\times \mathcal{U}_{t,T},$
$$u(t,\z)\geq
u(t,\eta)-\frac{\e}{3},\mbox{ \ }Y^{t,\eta;v}_{t}-Y^{t,\z;v}_{t}\geq
-\frac{\e}{3}, P-a.s..\eqno{(2.7)}$$

On the other hand, there exists $v^{i}(\cdot)\in \mathcal{U}_{t,T}$,
such that,
$$
u(t,x_{i})\geq Y^{t,x_{i};v^{i}.}-\frac{\e}{3}, i=1,2,...N_{\d}.
$$
So
$$
u(t,\eta)=\sum_{i=1}^{N_{\d}}I_{A_{i}}u(t,x_{i})\geq
\sum_{i=1}^{N_{\d}}I_{A_{i}}Y^{t,x_{i};v^{i}.}-\frac{\e}{3}=Y^{t,\eta;v.}-\frac{\e}{3},
$$
where $v(\cdot)=\sum_{i=1}^{N_{\d}}I_{A_{i}}v^{i}(\cdot).$ This with
(2.7), we have
$$
u(t,\z)\geq Y^{t,\z;v}_{t}-\e.
$$
$$
\eqno{\Box}
$$

Before stating the generalized Dynamic Programming Principle, let us
recall the following basic estimate of the solutions to BSDEs which
will be used often in the sequel(see Theorem 2.3. in [6]):

{\bf Lemma 2.4.} Consider the following two BSDEs:
$$
Y^{1}_{t}=\xi^{1}+\int_{t}^{T}[g(s,Y^{1}_{s},Z^{1}_{s})+\varphi^{1}_{s}]ds-\int_{t}^{T}Z^{1}_{s}dW_{s},\eqno{(a)}
$$
$$
Y^{2}_{t}=\xi^{2}+\int_{t}^{T}[g(s,Y^{2}_{s},Z^{2}_{s})+\varphi^{2}_{s}]ds-\int_{t}^{T}Z^{2}_{s}dW_{s},\eqno{(b)}
$$
where $\xi^{1},\xi^{2}\in L^{2}(\Omega,\mathcal{F}_{T},P;R^{m})$,
$\varphi^{1},\varphi^{2}\in \mathcal{M}(t,T;R^{m})$, and $g :
\Omega\times [0,T]\times R^{m}\times R^{m\times d}\rightarrow R^{m}$
satisfies: $\forall (y,z)\in R^{m}\times R^{m\times d}$,
$g(\cdot,y,z)$ is a $\mathcal{F}_{t}-$adapted process valued in
$R^{m}$ and
$$
\int_{0}^{T}|g(\cdot,0,0)|ds\in
L^{2}(\Omega,\mathcal{F}_{T},P;R^{m}),
$$
$$
|g(t,y,z)-g(t,y',z')|\leq C_{L}(|y-y'|+|z-z'|).
$$
Then the difference between the solutions to BSDE (a) and (b)
satisfies:
$$
\begin{array}{ll}
&|Y^{1}_{t}-Y^{2}_{t}|^{2}+\frac{1}{2}E^{\mathcal{F}_{t}}\int_{t}^{T}[|Y^{1}_{s}-Y^{2}_{s}|^{2}+|Z^{1}_{s}-Z^{2}_{s}|^{2}]e^{\beta_{0}(s-t)}ds\\
\leq&
E^{\mathcal{F}_{t}}|\xi^{1}-\xi^{2}|^{2}e^{\beta_{0}(T-t)}+E^{\mathcal{F}_{t}}\int_{t}^{T}|\varphi^{1}_{s}-\varphi^{2}_{s}|^{2}e^{\beta_{0}(s-t)}ds,
\end{array}
$$
where $\beta_{0}=16(1+C_{L}^{2})$.

Now let's consider the so-called backward semigroup (see [6]):
$\forall (t,x)\in [0,T]\times M, 0\leq\d\leq T-t, \eta\in
L^{2}(\Omega, \mathcal{F}_{t+\d},P;R)$, we set
$$
G^{t,x;v.}_{t,t+\d}[\eta]:=Y_{t},
$$
where $(Y_{s},Z_{s})_{t\leq s\leq t+\d}$ is the unique solution to
the following BSDE:
$$
\left\{\begin{array}{l}
-dY_{s}=f(s,X_{s}^{t,x;v.},Y_{s},Z_{s},v_{s})-Z_{s}dW_{s},s\in
[t,t+\d],
\\
Y_{t+\d}=\eta.
\end{array}\right.
$$
So obviously,
$$
G_{t,T}^{t,x;v.}[\P(X_{T}^{t,x;v.})]=G_{t,t+\d}^{t,x;v.}[Y^{t,x;v.}_{t+\d}].
$$
The follows is the generalized Dynamic Programming Principle (DPP in
short):

{\bf Theorem 2.5. }$\forall (t,x)\in [0,T]\times M,\forall \d\in
[0,T-t]$, we have
$$
\begin{array}{ll}
u(t,x)&=\mbox{essinf}_{v(\cdot)\in
\mathcal{U}_{t,t+\d}}G_{t,t+\d}^{t,x;v.}[u(t+\d,X^{t,x;v.}_{t+\d})]\\
&=\inf_{v(\cdot)\in\mathcal{U}_{t,t+\d}^{t}}G_{t,t+\d}^{t,x;v.}[u(t+\d,X^{t,x;v.}_{t+\d})].
\end{array}
\eqno{(2.8)}
$$

{\bf Proof: }We will only prove the first equality. By the
definition of $u(t,x)$, we have
$$
\begin{array}{ll}
u(t,x)&=\mbox{essinf}_{v(\cdot)\in
\mathcal{U}_{t,T}}G_{t,T}^{t,x;v.}[\P(X_{T}^{t,x;v.})]\\
&=\mbox{essinf}_{v(\cdot)\in
\mathcal{U}_{t,T}}G_{t,t+\d}^{t,x;v.}[Y_{t+\d}^{t+\d,X_{t+\d}^{t,x;v.};v.}].
\end{array}
$$
So by Lemma 2.3. and the comparison theorem of BSDEs, we get
$$
u(t,x)\geq \mbox{essinf}_{v(\cdot)\in
\mathcal{U}_{t,t+\d}}G_{t,t+\d}^{t,x;v.}[u(t+\d,X^{t,x;v.}_{t+\d})].
$$

On the other hand, for any $\e>0$, there exists an admissible
control $\bar{v}(\cdot)\in \mathcal{U}_{t+\d,T}$, such that,
$$
u(t+\d,X^{t,x;v.}_{t+\d})\geq
Y_{t+\d}^{t+\d,X_{t+\d}^{t,x;v.};\bar{v}.}-\e.
$$
Also by the comparison theorem of BSDEs and Lemma 2.4., we have
$$
\begin{array}{ll}
u(t,x)&=\mbox{essinf}_{v(\cdot)\in
\mathcal{U}_{t,T}}G_{t,t+\d}^{t,x;v.}[Y_{t+\d}^{t+\d,X_{t+\d}^{t,x;v.};v.}]\\
&\leq \mbox{essinf}_{v(\cdot)\in
\mathcal{U}_{t,t+\d}}G_{t,t+\d}^{t,x;v.}[Y_{t+\d}^{t+\d,X_{t+\d}^{t,x;v.};\bar{v}.}]\\
&\leq \mbox{essinf}_{v(\cdot)\in
\mathcal{U}_{t,t+\d}}G_{t,t+\d}^{t,x;v.}[u(t+\d,X^{t,x;v.}_{t+\d})+\e]\\
&\leq \mbox{essinf}_{v(\cdot)\in
\mathcal{U}_{t,t+\d}}G_{t,t+\d}^{t,x;v.}[u(t+\d,X^{t,x;v.}_{t+\d})]+C\e.
\end{array}
$$
Since $\e$ is arbitrary, (2.8) holds true.
$$
\eqno{\Box}
$$

In Lemma 2.2., we know that $u(t,x)$ is uniformly continuous in $x$,
uniformly in $t$. Now with the DPP, we can also get that $u(t,x)$ is
continuous in $t$.

{\bf Proposition 2.6.} The value function $u(t,x)$ is continuous in
$t,t\in [0,T]$.

{\bf Proof: }$\forall (t,x)\in [0,T]\times M, \d\in [0,T-t]$, by the
DPP, we have: $\forall \e>0$, there exists an admissible control
$\bar{v}(\cdot)\in \mathcal{U}^{t}_{t,T}$, such that,
$$
G^{t,x;\bar{v}.}_{t,t+\d}[u(t+\d,X^{t,x;\bar{v}.}_{t+\d})]\geq
u(t,x)\geq
G^{t,x;\bar{v}.}_{t,t+\d}[u(t+\d,X^{t,x;\bar{v}.}_{t+\d})]-\e.
\eqno{(2.9)}
$$
So
$$
u(t,x)-u(t+\d,x)\leq
G^{t,x;\bar{v}.}_{t,t+\d}[u(t+\d,X^{t,x;\bar{v}.}_{t+\d})]-u(t+\d,x)=I^{1}_{\d}+I^{2}_{\d},
\eqno{(2.10)}
$$
where
$$
\begin{array}{l}
I^{1}_{\d}=I^{1}_{\d}(\bar{v}(\cdot))=G^{t,x;\bar{v}.}_{t,t+\d}[u(t+\d,X^{t,x;\bar{v}.}_{t+\d})]-G^{t,x;\bar{v}.}_{t,t+\d}[u(t+\d,x)],\\
I^{2}_{\d}=I^{2}_{\d}(\bar{v}(\cdot))=G^{t,x;\bar{v}.}_{t,t+\d}[u(t+\d,x)]-u(t+\d,x).
\end{array}
$$

Now let us evaluate $I^{1}_{\d}$ and $I^{2}_{\d}$. Still by Lemma
2.4., we have
$$
|I^{1}_{\d}|\leq
[C_{0}E|u(t+\d,X^{t,x;\bar{v}.}_{t+\d})-u(t+\d,x)|^{2}]^{\frac{1}{2}}.
\eqno{(2.11)}
$$
We can use the similar method as in (3.6) to get that
$$
E[\sup_{v\in
\mathcal{U}^{t}_{t,T}}d^{2}(X^{t,x;v.}_{t+\d},x)]\rightarrow 0,
\mbox{ \ when \ }\d\rightarrow 0.
$$
Thus for all $\e_{0}>0$,
$$
\lim_{\d\rightarrow 0}\sup_{v\in
\mathcal{U}^{t}_{t,T}}P\{d^{2}(X^{t,x;v.}_{t+\d},x)>\e_{0}\}=0.
$$
By the continuity of $u$ w.r.t $x$, we have
$$
\lim_{\d\rightarrow 0}\sup_{v\in \mathcal{U}^{t}_{t,T}}\sup_{s\in
[t,t+\d]}P\{|u(s,X^{t,x;v.}_{t+\d})-u(s,x)|^{2}>\e_{0}\}=0.
$$
Recall that $u(t,x)$ is bounded, so
$$
\lim_{\d\rightarrow 0}\sup_{v\in \mathcal{U}^{t}_{t,T}}\sup_{s\in
[t,t+\d]}E|u(s,X^{t,x;v.}_{t+\d})-u(s,x)|^{2}=0.
$$
Through (2.11), we know that $$\lim_{\d\rightarrow
0}[\sup_{v\in\mathcal{U}^{t}_{t,T}}|I^{1}_{\d}(v(\cdot))|]=0.$$

For $I^{2}_{\d}$, we have
$$
\begin{array}{ll}
I^{2}_{\d}&=E[u(t+\d,x)+\int_{t}^{t+\d}f(s,X^{t,x;\bar{v}.}_{s},Y^{t,x;\bar{v}.}_{s},Z^{t,x;\bar{v}.}_{s},\bar{v}_{s})ds
-\int_{t}^{t+\d}Z^{t,x;\bar{v}.}_{s}dW_{s}]-u(t+\d,x)\\
&=E[\int_{t}^{t+\d}f(s,X^{t,x;\bar{v}.}_{s},Y^{t,x;\bar{v}.}_{s},Z^{t,x;\bar{v}.}_{s},\bar{v}_{s})ds].
\end{array}
$$
From the assumptions (A1) and (A2),
$$
\begin{array}{ll}
|I^{2}_{\d}|&\leq
E\int^{t+\d}_{t}K(|Y^{t,x;\bar{v}.}_{s}|+|Z^{t,x;\bar{v}.}_{s}|)ds+K_{0}\d\\
&\leq
C(\int_{t}^{t+\d}E(|Y^{t,x;\bar{v}.}_{s}|^{2}+|Z^{t,x;\bar{v}.}_{s}|^{2})e^{\beta_{1}(s-t)}ds)^{\frac{1}{2}}+K_{0}\d,
\end{array}
$$
where $\beta_{1}$ is a constant defined in Lemma 2.2. Note (2.6) in
Lemma 2.2., we have
$$
\lim_{\d\rightarrow
0}[\sup_{v\in\mathcal{U}^{t}_{t,T}}|I^{2}_{\d}(v(\cdot))|]=0.
$$
From (2.10), we can conclude that
$$
\limsup_{\d\rightarrow 0}[u(t,x)-u(t+\d,x)]\leq 0. \eqno{(2.12)}
$$

Consider the right inequality of (2.9), we have
$$
u(t,x)-u(t+\d,x)\geq
G^{t,x;\bar{v}.}_{t,t+\d}[u(t+\d,X^{t,x;\bar{v}.}_{t+\d})]-u(t+\d,x)-\e=I^{1}_{\d}+I^{2}_{\d}-\e.
$$
Thus
$$
\liminf_{\d\rightarrow 0}[u(t,x)-u(t+\d,x)]\geq
\liminf_{\d\rightarrow
0}[-\sup_{v\in\mathcal{U}^{t}_{t,T}}|I^{1}_{\d}(v(\cdot))|-\sup_{v\in\mathcal{U}^{t}_{t,T}}|I^{2}_{\d}(v(\cdot))|]-\e.
$$
So for any $\e>0$,
$$
\liminf_{\d\rightarrow 0}[u(t,x)-u(t+\d,x)]\geq -\e.
$$
That is
$$
\liminf_{\d\rightarrow 0}[u(t,x)-u(t+\d,x)]\geq 0.
$$
This with (2.12), we get that $\lim_{\d\rightarrow
0}[u(t,x)-u(t+\d,x)]=0$.
$$
\eqno{\Box}
$$

{\bf Remark 2.7.} We can conclude from this proposition and Lemma
2.2. that $u(t,x)$ is continuous in $(t,x)\in [0,T]\times M$.

\section{The viscosity solution to
the generalized Hamilton-Jacobi-Bellman equation on Rimeannian
manifolds}

As it is well known to all, the value function $u(t,x)$ is usually a
viscosity solution to some Hamilton-Jacobi-Bellman(H-J-B in short)
equation. Our generalized H-J-B equation is:
$$
\left\{
\begin{array}{l}
\partial_{t}u(t,x)+\inf_{v\in
U}\{(v_{0}V_{0}u)(t,x)+\frac{1}{2}\sum_{\a=1}^{d}(v_{\a}^{2}V_{\a}V_{\a}u)(t,x)+f(t,x,u,\{v_{\a}V_{\a}u\}_{\a=1}^{d},v)\}=0,\\
u(T,x)=\P(x),
\end{array}
 \right.
 \eqno{(3.1)}
$$
where $v=(v_{0},v_{1},...,v_{d})\in U$ and we denote by
$\{v_{\a}V_{\a}u\}_{\a=1}^{d}=(v_{1}V_{1}u,...,v_{d}V_{d}u)$ an
element in $R^{1\times d}$. This is a fully nonlinear second order
parabolic PDEs on Riemannian manifolds.

The theory of viscosity solutions to PDEs in Euclidean space was
introduced by M. G. Crandall and P. L. Lions in the 1980's (see
[2]). Until recently, D. Azagra, J. Ferrera and B. Sanz [1] gave a
work about Dirichlet problem on a complete Riemannian manifold with
some restrictions on curvature. X. Zhu [7] studied parabolic PDEs on
Riemannian manifolds.

{\bf Definition 3.1. }We say $u\in C([0,T]\times M)$ is a viscosity
supersolution (subsolution) of (3.1), if $u(T,x)\geq \P(x)$ ($\leq
\P(x)$), and for all $\varphi\in C^{1,2}([0,T]\times M)$, at each
minimum (maximum) point $(t,x)$ of $u-\varphi$ and
$u(t,x)=\varphi(t,x)$, the following inequality holds:
$$
\partial_{t}\varphi(t,x)+\inf_{v\in
U}\{(v_{0}V_{0}\varphi)(t,x)+\frac{1}{2}\sum_{\a=1}^{d}(v_{\a}^{2}V_{\a}V_{\a}\varphi)(t,x)+f(t,x,u,\{v_{\a}V_{\a}\varphi\}_{\a=1}^{d},v)\}\leq0
(\geq 0).
$$
$u(t,x)$ is said to be a viscosity solution of (3.1) if it is both a
viscosity supersolution and a subsoultion.

We have proved that $u(t,x)$ is continuous in $(t,x)\in [0,T]\times
M$. So we are ready to present that $u(t,x)$ is a viscosity solution
of (3.1).

{\bf Theorem 3.2. }Under the assumptions (A1) and (A2), the value
function $u(t,x)$ is a viscosity solution to H-J-B equation (3.1).

To prove this theorem, we need the following three lemmas. We set
$$
\begin{array}{ll}
F(s,x,y,z,v) =&\partial_{t}\varphi(s,x)+(v_{0}V_{0}\varphi)(s,x)
+\frac{1}{2}\sum_{\a=1}^{d}(v_{\a}^{2}V_{\a}V_{\a}\varphi)(s,x)\\
&+f(s,x,y+\varphi(s,x),z+\{v_{\a}V_{\a}\varphi\}_{\a=1}^{d},v),
\end{array}
$$
and consider the a BSDE defined on $[t,t+\d]$:
$$
\left\{\begin{array}{l} -dY_{s}^{1;v.}=F(s,X_{s}^{t,x ;v.},Y_{s}^{1
;v.},Z_{s}^{1 ;v.},v_{s})-Z_{s}^{1;v.}dW_{s},
\\
Y_{t+\d}^{1 ;v.}=0,
\end{array}\right.
\eqno{(3.2)}
$$
Then let us consider following first lemma:

{\bf Lemma 3.3. }$\forall s\in[t,t+\d]$, we have
$$
Y^{1;v.}_{s}=G_{s,t+\d}^{t,x,v}[\varphi(t+\d,X^{t,x;v.}_{t+\d})]-\varphi(s,X^{t,x;v.}_{s}),a.s..
$$

{\bf Proof: }Recall that
$G_{s,t+\d}^{t,x,v}[\varphi(t+\d,X^{t,x;v.}_{t+\d})]$ is defined
through the solution to the following BSDE:
$$
\left\{\begin{array}{l} -dY_{s}^{v.}=f(s,X_{s}^{t,x
;v.},Y_{s}^{v.},Z_{s}^{v.},v_{s})-Z_{s}^{v.}dW_{s},s\in[t,t+\d],
\\
Y_{t+\d}^{v.}=\varphi(t+\d,X^{t,x;v.}_{t+\d}).
\end{array}\right.
\eqno{(3.3)}
$$
That is
$$
G_{s,t+\d}^{t,x,v}[\varphi(t+\d,X^{t,x;v.}_{t+\d})]=Y^{v.}_{s}, s\in
[t,t+\d].
$$
So what we need to do is proving that
$Y_{s}^{1;v.}+\varphi(s,X_{s}^{t,x;v.}),s\in [t,t+\d]$ is also a
solution to (3.3).

Applying It\^{o}'s formula to
$Y_{s}^{1;v.}+\varphi(s,X_{s}^{t,x;v.})$:
$$
\begin{array}{ll}
& -d(Y_{s}^{1;v.}+\varphi(s,X_{s}^{t,x;v.}))\\
=&
F(s,X_{s}^{t,x;v.},Y^{1;v}_{s},Z^{1;v}_{s},v_{s})ds-Z^{1;v}_{s}dW_{s}\\
&-[\partial_{t}\varphi(s,X_{s}^{t,x;v.})+(v_{0}V_{0}\varphi)(s,X_{s}^{t,x;v.})+\frac{1}{2}\sum_{\a=1}^{d}(v_{\a}^{2}V_{\a}V_{\a}\varphi)(s,X_{s}^{t,x;v.})]ds\\
&-\{(v_{\a}V_{\a}\varphi)(s,X_{s}^{t,x;v.})\}_{\a=1}^{d}dW_{s}\\
=&f(s,X_{s}^{t,x;v.},Y_{s}^{1;v.}+\varphi(s,X_{s}^{t,x;v.}),Z^{1;v}_{s}+\{(v_{\a}V_{\a}\varphi)(s,X_{s}^{t,x;v.})\}_{\a=1}^{d})\\
&-(Z^{1;v}_{s}+\{(v_{\a}V_{\a}\varphi)(s,X_{s}^{t,x;v.})\}_{\a=1}^{d})dW_{s}.
\end{array}
$$
Moreover,
$$
(Y_{s}^{1;v.}+\varphi(s,X_{s}^{t,x;v.}))_{s=t+\d}=\varphi(t+\d,X_{t+\d}^{t,x;v.}).
$$
So $Y_{s}^{1;v.}+\varphi(s,X_{s}^{t,x;v.}),s\in [t,t+\d]$ is in fact
a solution to (3.3). By the uniqueness of the solution to (3.2),
we've finished the proof.
$$
\eqno{\Box}
$$

Consider the following BSDE which is easier than (3.2):
$$
\left\{\begin{array}{l} -dY_{s}^{2;v.}=F(s,x,Y_{s}^{2 ;v.},Z_{s}^{2
;v.},v_{s})-Z_{s}^{2;v.}dW_{s},s\in [t,t+\d],
\\
Y_{t+\d}^{2 ;v.}=0.
\end{array}\right.
\eqno{(3.4)}
$$
The following lemma shows that, when $\d$ is small enough, the
difference between BSDE(3.2) and (3.4) $|Y_{t}^{1;v.}-Y_{t}^{2;v.}|$
can be ignored.

{\bf Lemma 3.4.} We have the following estimate:
$$
|Y^{1;v.}_{t}-Y^{2;v.}_t|\leq C\d\rho_{1}(\d),\eqno{(3.5)}
$$
where, $\rho_{1}(\d)\downarrow 0$, when $\d\downarrow 0$, and
$\rho_{1}(\cdot)$ does not depend on the control $v(\cdot)\in
\mathcal{U}^{t}_{t,t+\d}$.

{\bf Proof: }If we set
$$\eta^{\d}:=\sup_{s\in[t,t+\d]}d(X^{t,x;v.}_{s},x).$$ With the
similar method as in (3.6), We have
$$
E[\eta^{\d}]\downarrow 0, \mbox{ \ when \ }\d\downarrow
0.\eqno{(3.6)}
$$

We can use the estimate in Lemma 2.4. to BSDE(3.2) and (3.4) with
$\xi^{1}=\xi^{2}=0$,
$$
g(s,y,z)=F(s,X^{t,x;v.}_{s},y,z,v_{s}),
$$
$$
\varphi^{1}_{s}=0,\varphi^{2}_{s}=F(s,x,Y^{2;v.}_{s},Z^{2;v.}_{s},v_{s})-F(s,X^{t,x;v.}_{s},Y^{2;v.}_{s},Z^{2;v.}_{s},v_{s}).
$$
From the definition of the function $F$, we can see that $g(s,y,z)$
satisfies Lipschitz conditions w.r.t. $(y,z)$. And there exists a
constant $C$ and a function $\rho(\cdot)$ with
$$\rho(\d)\downarrow 0,\mbox{ \ when \ }\d\downarrow 0$$ which do not
depend on $v(\cdot)$, such that,
$$
|\varphi^{2}_{s}|\leq C\rho(d(X^{t,x;v.}_{s},x)),
$$
where, when $\e\downarrow 0$, $\rho(\e)\rightarrow 0$.

According Lemma 2.4., we have
$$
E\int_{t}^{t+\d}[|Y^{1;v.}_{s}-Y^{2;v.}_{s}|^{2}+|Z^{1;v.}_{s}-Z^{2;v.}_{s}|^{2}]ds
\leq CE\int_{t}^{t+\d}\rho^{2}(d(X^{t,x;v.}_{s},x))ds\leq C\d
E\rho^{2}(\eta^{\d}).\eqno{(3.7)}
$$
On the other hand, since $Y^{1;v.}_{t}$ and $Y^{2;v.}_{t}$ are both
deterministic when $v(\cdot)\in \mathcal{U}^{t}_{t,t+\d}$, apply
It\^{o}'s  formula to $Y^{1;v.}_{s}-Y^{2;v.}_{s}$ on $[t,t+\d]$, we
have
$$
\begin{array}{ll}
& |Y^{1;v.}_{t}-Y^{2;v.}_{t}|=|E(Y^{1;v.}_{t}-Y^{2;v.}_{t})|\\
\leq
&E\int_{t}^{t+\d}|F(s,X^{t,x;v.}_{s},Y^{1;v.}_{s},Z^{1;v.}_{s},v_{s})-F(s,x,Y^{2;v.}_{s},Z^{2;v.}_{s},v_{s})|ds\\
\leq &
E\int_{t}^{t+\d}[K|Y^{1;v.}_{s}-Y^{2;v.}_{s}|+K|Z^{1;v.}_{s}-Z^{2;v.}_{s}|+C\rho(d(X^{t,x;v.}_{s},x))]ds\\
\leq & C\d
E\rho(\eta^{\d})+C\d^{\frac{1}{2}}\{E\int_{t}^{t+\d}[|Y^{1;v.}_{s}-Y^{2;v.}_{s}|^{2}+|Z^{1;v.}_{s}-Z^{2;v.}_{s}|^{2}]ds\}^{\frac{1}{2}}.
\end{array}
$$
This with (3.7), we get
$$
|Y^{1;v.}_{t}-Y^{2;v.}_{t}|\leq
C\d[E\rho(\eta^{\d})+\{E\rho^{2}(\eta^{\d})\}^{\frac{1}{2}}].
$$
Because of the compactness of $M$, there exists a constant $C$, such
that, $ \eta^{\d}\leq C,\forall \d\geq 0$. Thus
$$E\rho^{2}(\eta^{\d})<\infty.$$ This, together with (3.6) we have
that (3.5) holds true.
$$
\eqno{\Box}
$$

The following lemma tells us how to compute $\inf_{v(\cdot)\in
\mathcal{U}^{t}_{t,t+\d}}Y^{2;v.}_{t}$:

{\bf Lemma 3.5. }We have
$$
\inf_{v(\cdot)\in \mathcal{U}^{t}_{t,t+\d}}Y^{2;v.}_{t}=Y_{0}(t),
$$
where $Y_{0}(t)$ is the solution to the following ODE:
$$
\left\{
\begin{array}{l}
-\dot{Y}_{0}(s)=F_{0}(s,x,Y_{0}(s),0),s\in [t,t+\d],\\
Y_{0}(t+\d)=0,
\end{array}
\right.\eqno{(3.8)}
$$
and the function $F_{0}$ is defined as:
$$
F_{0}(t,x,y,z):=\inf_{v\in U}F(t,x,y,z,v).
$$

{\bf Proof: }Consider the following BSDE:
$$
\left\{
\begin{array}{l}
-dY^{0}_{s}=F_{0}(s,x,Y^{0}_{s},Z^{0}_{s})ds-Z^{0}_{s}dW_{s},s\in [t,t+\d],\\
Y^{0}_{t+\d}=0.
\end{array}
\right.\eqno{(3.9)}
$$
Note that $F_{0}$ is a deterministic function of $(t,x,y,z)$. So the
solution to (3.9) is just:
$$
(Y^{0}_{s},Z^{0}_{s})=(Y_{0}(s),0),s\in [t,t+\d],
$$
that is to say equation (3.8) and BSDE (3.9) are the same one.

By the definition of $F_{0}$, we know
$$
F_{0}(s,x,y,z)\leq F(s,x,y,z,v_{s}),\mbox{ \ \ }\forall v(\cdot)\in
\mathcal{U}^{t}_{t,t+\d},s\in [t,t+\d].
$$
Through the comparison theorem of the solutions to BSDE (3.4) and
(3.9), $\forall v(\cdot)\in \mathcal{U}^{t}_{t,t+\d}$,
$$
Y^{0}_{s}\leq Y^{2;v.}_{s},s\in [t,t+\d].
$$
So
$$
Y_{0}(t)=Y^{0}_{t}\leq \inf_{v(\cdot)\in
\mathcal{U}^{t}_{t,t+\d}}Y^{2;v.}_{t}.
$$

On the other hand, if we denote by $\mathcal{U}^{0}_{t,t+\d}$ the
set of all admissible controls in $[t,t+\d]$ which are deterministic
processes. Then we can show that
$$
Y^{0}_{t}\geq \inf_{v(\cdot)\in
\mathcal{U}^{0}_{t,t+\d}}Y^{2;v.}_{t}.
$$
That is because, $\forall v(\cdot)\in \mathcal{U}^{0}_{t,t+\d}$,
$(Y^{2;v.}_{s},s\in [t,t+\d])$ is the solution to the following ODE:
$$
\left\{
\begin{array}{l}
-\dot{Y}^{2;v.}(s)=F(s,x,Y^{2;v.}(s),0,v_{s}),s\in [t,t+\d],\\
Y^{2;v.}(t+\d)=0.
\end{array}
\right.
$$
According to the definition of $F$, for all $\e>0$, there exists
$(v_{s}(\e),s\in [t,t+\d])\in \mathcal{U}^{0}_{t,t+\d}$, such that,
$$
F_{0}(s,x,y,z)\geq F(s,x,y,z,v_{s}(\e))-\frac{\e}{\d},s\in [t,t+\d].
$$
That yields
$$
Y^{0}_{t}\geq Y^{2;v.(\e)}_{t}-\e.
$$
Since $\e$ is arbitrary, we have
$$
Y^{0}_{t}\geq \inf_{v(\cdot)\in
\mathcal{U}^{0}_{t,t+\d}}Y^{2;v.}_{t}.
$$
And consequently
$$
Y_{0}(t)=Y^{0}_{t}\geq\inf_{v(\cdot)\in
\mathcal{U}^{0}_{t,t+\d}}Y^{2;v.}_{t}\geq \inf_{v(\cdot)\in
\mathcal{U}^{t}_{t,t+\d}}Y^{2;v.}_{t}.
$$
$$
\eqno{\Box}
$$

After the above preparation, we can show the proof of Theorem 3.2.

{\bf The proof of Theorem 3.2.}  For all $\varphi\in
C^{1,2}([0,T\times M])$, suppose that $(t,x)$ is a minimum(resp.,
maximum) point of $u-\varphi$ and $u(t,x)-\varphi(t,x)=0$. By
DPP(2.8), we have
$$
\varphi(t,x)=u(t,x)=\inf_{v(\cdot)\in
\mathcal{U}^{t}_{t,T}}G^{t,x;v.}_{t,t+\d}[u(t+\d,X^{t,x;v.}_{t+\d})].
$$
Since
$$
u(t+\d,X^{t,x;v.}_{t+\d})\geq \varphi(t+\d,X^{t,x;v.}_{t+\d}),
$$
this with the comparison theorem of the solutions to BSDEs,
$$
\inf_{v(\cdot)\in
\mathcal{U}^{t}_{t,T}}\{G^{t,x;v.}_{t,t+\d}[\varphi(t+\d,X^{t,x;v.}_{t+\d})]-\varphi(t,x)\}\leq
0\mbox{ \ (resp., }\geq 0).
$$
Through Lemma 3.3., we have
$$
\inf_{v(\cdot)\in \mathcal{U}^{t}_{t,T}}Y^{1;v.}_{t}\leq 0\mbox{ \
(resp., }\geq 0).
$$
So by (3.5),
$$
\inf_{v(\cdot)\in \mathcal{U}^{t}_{t,T}}Y^{2;v.}_{t}\leq
C\d\rho_{1}(\d)\mbox{ \ (resp., }\geq -C\d\rho_{1}(\d)).
$$
According Lemma 3.5., it yields
$$
Y_{0}(t)\leq C\d\rho_{1}(\d)\mbox{ \ (resp., }\geq
-C\d\rho_{1}(\d)).
$$
So
$$
\begin{array}{ll}
 Y_{0}(t)&=\int_{t}^{t+\d}F_{0}(s,x,Y_{0}(s),0)ds\\
 &=\d F_{0}(t+\d,x,Y_{0}(t+\d),0)+o(\d)\\
 &=\d F_{0}(t+\d,x,0,0)+o(\d)\\
&\leq C\d\rho_{1}(\d)\mbox{ \ (resp., }\geq -C\d\rho_{1}(\d)).
\end{array}
$$
Divided by $\d$ and let $\d\downarrow 0$, we have
$$
F_{0}(t,x,0,0)=\inf_{v\in U}F(t,x,0,0,v)\leq 0\mbox{ \ (resp., }\geq
0).
$$
That is to say
$$
\begin{array}{r}
\partial_{t}\varphi(t,x)+\inf_{v\in
U}\{(v_{0}V_{0}\varphi)(t,x)+\frac{1}{2}\sum_{\a=1}^{d}(v_{\a}^{2}V_{\a}V_{\a}\varphi)(t,x)+f(t,x,u,\{v_{\a}V_{\a}\varphi\}_{\a=1}^{d},v)\}\leq0,\\
(\mbox{resp,. }\geq 0).
\end{array}
$$ And obviously, $u(T,x)=\P(x)$. So $u(t,x)$ is a viscosity
solution to PDE (3.1).

Now let's deal with the uniqueness conclusion. We assume that for
all $x,y\in M$ s.t. $d(x,y)<\min\{i_{M}(x),i_{M}(y)\}$, $t\in
[0,T]$,
$$
\begin{array}{ll}
(H1) & \|L_{xy}V_{0}(t,x)-V_{0}(t,y)\| \leq \mu
d(x,y),\\
(H2) & L_{xy}V_{\a}(t,x)= V_{\a}(t,y),\a=1,2,...,d,
\end{array}
$$
where $\mu$ is a positive constant.

Consider a generalized case:
$$\left\{
\begin{array}{l}
u_{t}+\inf_{v\in U}H(t,x,u,du,d^{2}u,v)=0\mbox{ \ in \ }(0,T)\times M,\\
u(0,x)=\psi(x),  x\in M,
\end{array}
\right. \eqno{(3.10)}
$$
where, $du, d^{2}u$ mean $d_{x}u(t,x)$ and $d^{2}_{x}u(t,x)$.

Set
$$
\chi:=\{(t,x,r,\z,A,v): t\in [0,T], x\in M, r\in R, \z\in
TM^{*}_{x}, A\in {\cal{L}}^{2}_{s}(TM_{x}),v\in U\},
$$
where $TM_{x}^{*}$ stands for the cotangent space of $M$ at a point
$x$, $TM_{x}$ stands for the tangent space at $x$ and
${\cal{L}}_{s}^{2}(TM_{x})$ denotes the symmetric bilinear forms on
$TM_{x}$.

From [7], we have the following comparison theorem of viscosity
solutions to PDE (3.10):

 {\bf Theorem 3.6. }Let $M$ be a compact Riemannian manifold
(without boundary) , and $H : \chi \rightarrow R$ be continuous,
proper for each fixed $(t,v)\in (0,T)\times U$ and satisfy: there
exists a function $\omega : [0,+\infty]\rightarrow [0,+\infty]$ with
$\omega(0+)=0$ and such that
$$
\sup_{v\in
U}[H(t,y,r,\a\exp_{y}^{-1}(x),Q,v)-H(t,x,r,-\a\exp_{x}^{-1}(y),P,v)]\leq
\omega(\a d^{2}(x,y)+d(x,y)),\eqno{(3.11)}
$$
for all fixed $t\in (0,T)$ and for all $x, y\in M, r\in R, P\in
T_{2,s}(M)_{x}, Q\in T_{2,s}(M)_{y}$ with
$$
-(\frac{1}{\e_{\a}}+\|A_{\a}\|)\left(
\begin{array}{cc}
I & 0\\
0 & I
\end{array}
\right)\leq \left(
\begin{array}{cc}
P & 0\\
0 & -Q
\end{array}
\right)\leq A_{\a}+\e_{\a}A_{\a}^{2}, \eqno{(3.12)}
$$
where $A_{\a}$ is the second derivative of the function
$\p_{\a}(x,y)=\frac{\a}{2}d^{2}(x,y)(\a >0)$ at the point $(x,y)\in
M\times M$,
$$
\e_{\a}=\frac{1}{2(1+\|A_{\a}\|)}
$$
and the points $x, y$ are assumed to be close enough to each other
so that  $d(x,y)<\min\{i_{M}(x),i_{M}(y)\}$.

Let $u\in USC([0,T)\times M)$ be a subsolution and $v\in
LSC([0,T)\times M)$ a supersolution of (3.10). Then $u\leq v$ on
$[0,T)\times M$. In particular PDEs (3.10) has at most one viscosity
solution.

{\bf Theorem 3.7.} The value function $u(t,x)$ is the unique
viscosity solution to PDE(3.1).

{\bf Proof: } For PDE(3.1), the function $H : [0,T]\times M\times
R\times TM^{*}_{x}\times {\cal{L}}^{2}_{s}(TM_{x})\times
U\rightarrow R$ is:
$$
\begin{array}{ll}
&H(t,x,r,\zeta, P,v)\\
=&-f(t,x,r,\langle \zeta,
v_{\a}V_{\a}(t,x)\rangle^{d}_{\a=1},v)-\langle \zeta,
v_{0}V_{0}(t,x)\rangle
-\frac{1}{2}\sum\limits_{\a=1}^{d}v^{2}_{\a}\langle
PV_{\a}(t,x),V_{\a}(t,x)\rangle.
\end{array}
$$
So for any fixed $t\in (0,T), r\in R$, when $d(x,y)
<\frac{1}{2}i_{M}$, if $P\in {\cal{L}}_{s}^{2}(TM_{x}), Q\in
{\cal{L}}_{s}^{2}(TM_{y})$ satisfy the following matrix inequality:
$$
-(\frac{1}{\e_{\a}}+\|A_{\a}\|)\left(
\begin{array}{cc}
I & 0\\
0 & I
\end{array}
\right)\leq \left(
\begin{array}{cc}
P & 0\\
0 & -Q
\end{array}
\right)\leq A_{\a}+\e_{\a}A_{\a}^{2},
$$
where $A_{\a}$ is the second derivative of the function
$\p_{\a}(x,y)=\frac{\a}{2}d^{2}(x,y)$ at the point $(x,y)\in M\times
M$,
$$
\e_{\a}=\frac{1}{2(1+\|A_{\a}\|)}.
$$
Since $M$ is compact, there is a $k_{0}>0$ s.t. the sectional
curvature is bounded below by $-k_{0}$ on $M$. So by Remark 4.7 in
[1],
$$ P-L_{yx}(Q)\leq \frac{3}{2}k_{0}\a d^{2}(x,y)I.$$
Then
$$
\begin{array}{rl}
& \sup\limits_{v\in U,t\in (0,T)}[H(t,y,r,\a\exp^{-1}_{y}(x),Q,v)-H(t,x,r,-\a\exp^{-1}_{x}(y),P,v)]\\
=&
 \sup\limits_{v\in U,t\in (0,T)}\{f(t,x,r,\langle
 -\a\exp^{-1}_{x}(y),v_{\a}V_{\a}(t,x)\rangle^{d}_{\a=1},v)\\
 & \ \ \ \ \ \ \ \  \ \ \ \ \ \ -f(t,y,r,\langle
 \a\exp^{-1}_{y}(x),v_{\a}V_{\a}(t,y)\rangle^{d}_{\a=1},v)\\
 & \ \ \ \ \ \ \ \  \ \ \ \ \ \  +\frac{1}{2}\sum\limits_{\a=1}^{d}v^{2}_{\a}[\langle
PV_{\a}(t,x),V_{\a}(t,x)\rangle-\langle
QV_{\a}(t,y),V_{\a}(t,y)\rangle]\\
 & \ \ \ \ \ \ \ \  \ \ \ \ \ \ -\langle \a\exp^{-1}_{x}(y),v_{0}V_{0}(t,x)\rangle
 -\langle\a\exp^{-1}_{y}(x),v_{0}V_{0}(t,y)\rangle\}\\
 \leq & \sup\limits_{v\in U,t\in (0,T)}\{Kd(x,y)
 +\frac{3}{4}k_{0}\a d^{2}(x,y)\sum\limits_{\a=1}^{d}v^{2}_{\a}\langle V_{\a},V_{\a}\rangle_{(t,x)}\\
& \ \ \ \ \ \ \ \  \ \ \ \ \ \  +\langle
L_{yx}\a\exp^{-1}_{y}(x),v_{0}V_{0}(t,x)\rangle
 -\langle\a\exp^{-1}_{y}(x),v_{0}V_{0}(t,y)\rangle\}\\
\leq &\bar{C}(\a d^{2}(x,y)+d(x,y)).
\end{array}
$$
Since $M$ and $U$ are both compact, through (A1), (H1) and (H2), we
get the last inequality. Where $\bar{C}$ is a constant positive. By
Theorem 3.6. we can get the uniqueness result of the viscosity
solution to PDE (3.1).
$$
\eqno{\Box}
$$


\begin{thebibliography}{00}

\bibitem{}D. Azagra, J. Ferrera and
B. Sanz, Viscosity solutions to second order partial differential
equations on Riemannian manifolds, \emph{Journal of Differential
Equations, Vol. 245(2008), 307-336.}

\bibitem{} M.G. Crandall, H. Ishii, P.-L. Lions, User's guide to
viscosity solutions of second order partial differential equations,
\emph{Bull, Amer. Math. Soc. 27 (1992) no. 1, 1-67.}

\bibitem{}N. El Karoui, S. Peng and M. C. Quenez, Backward
stochastic differential equation in finance, \emph{Math. Finance, 7
(1997), pp. 1-71}


\bibitem{} Elton P. Hsu, Stochastic Analysis on Manifolds,\emph{ American
Mathematical Society, 2002.}

\bibitem{} S. Peng, A generalized dynamic programming
principle and Hamilton-Jacobi-Bellmen equation, \emph{Stochastics
Stochastics Rep., 38(1992), 119-134.}

\bibitem{} J. Yan, S. Peng, S. Fang and L. Wu, Topics on
stochastic analysis, \emph{Science Press. Beijing (in Chinese),
1997.}

\bibitem{} X. Zhu, Viscosity solutions to second order parabolic
PDEs on Riemannian manifolds. Preprint.



\end{thebibliography}
\end{document}